\documentclass[a4paper,reqno]{amsart}
\usepackage{amssymb,amsmath,setspace,centernot}
\usepackage{makecell} %enable linebreak in tabular-environment
\newcommand{\Z}{\mathbb{Z}}
\newcommand{\Q}{\mathbb{Q}}
\newcommand{\R}{\mathbb{R}}
\newcommand{\C}{\mathbb{C}}
\newcommand{\N}{\mathbb{N}}

\newcommand\ee{\mathcal{E}}
\newcommand\Ce{\mathcal{C}}
\newcommand\sign{\operatorname{sgn}}
\newcommand\im{\operatorname{Im}}

\newcommand\smod[1]{\operatorname{mod} #1}
\newcommand\trans{\mathrm{t}}
\newcommand\aut[1]{\operatorname{Aut}^+(Q,\Z^{#1})}
\renewcommand\theta{\vartheta}
\newcommand{\modtheta}{\skew{0}\widehat{\vphantom{\rule{1pt}{8.4pt}}\smash{\widehat{\Theta}}}} %doublehat theta
\newcommand\ko{\hspace{.1em}}
\newcommand\steilt{\hspace{1pt}|\hspace{1pt}}
\newcommand\nsteilt{\hspace{1pt}\centernot{|}\hspace{1pt}}
\newtheorem{theorem}{Theorem}[section]
\newtheorem{lemma}[theorem]{Lemma}

\newtheorem{corollary}[theorem]{Corollary}
\theoremstyle{definition}
\newtheorem{definition}[theorem]{Definition}
\newtheorem{remark}[theorem]{Remark} 
\newtheorem{example}[theorem]{Example}
\newenvironment{proofof}[1]{\begin{proof}[Proof of {#1}]}{\end{proof}}

\makeatletter
\@namedef{subjclassname@2020}{\textup{2020} Mathematics Subject Classification}
\makeatother

\numberwithin{equation}{section}

\title[Indefinite Theta Series with (Spherical) Polynomials]{Theta Series for Quadratic Forms of Signature $(n-1,1)$\\ with (Spherical) Polynomials}
\date{\today}
\author[C. Roehrig, S. Zwegers]{Christina Roehrig and Sander Zwegers}
\address{Department of Mathematics and Computer Science, University of Cologne, Weyertal 86--90, 50931 Cologne, Germany}
\email{croehrig@math.uni-koeln.de}
\email{szwegers@uni-koeln.de}

\begin{document}

\begin{abstract}
We construct almost holomorphic and holomorphic modular forms by considering theta series for quadratic forms of signature $(n-1,1)$.
We include homogeneous and spherical polynomials in the definition of the theta series (generalizing a construction of the second author) to obtain holomorphic, almost holomorphic and modular theta series.
We give a criterion for these series to coincide, enabling us to %quickly
construct almost holomorphic and holomorphic cusp forms %of half-integral weight
on congruence subgroups of the modular group.
Further, we provide numerous explicit examples.
\end{abstract}

\maketitle

\section{Introduction}
One of the few known general constructions of holomorphic modular forms is via theta series:
if $Q:\R^n\longrightarrow \R$ is a positive definite quadratic form, which is integer-valued on the lattice $\Z^n$, and $P:\R^n\longrightarrow \C$ a spherical polynomial of degree $d$, then the theta series
\[\Theta_{Q,P} (\tau) := \sum_{\ell\in \Z^n} P(\ell)\ko q^{Q(\ell)}\qquad (q= e^{2\pi i\tau},\ \im(\tau)> 0)\]
is a (holomorphic) modular form of weight $n/2+d$ (on some subgroup of $\operatorname{SL}_2(\Z)$, with some character; see \cite{schoeneberg,ogg} for the case that $n$ is even and \cite{shimura} for the case that $n$ is odd).
If for example we take $n=4$, $Q=\|\cdot \|^2$ and $P\equiv 1$, then the theta function $\Theta_{Q,P}$ is modular of weight 2 on $\Gamma_0(4)$.
Writing this theta series as a linear combination of Eisenstein series we obtain Jacobi's four-square theorem, which gives a formula for the number of ways that a given positive integer $m$ can be represented as the sum of four squares:
\[ |\{ \ell\in \mathbb{Z}^4 \mid \ell_1^2+\ell_2^2+\ell_3^2+\ell_4^2 = m \}|= 8 \sum_{d \steilt m,\, 4 \nsteilt d} d\]

For indefinite quadratic forms the situation is more complicated, because the sum $\sum_{\ell\in\Z^n} P(\ell)\ko q^{Q(\ell)}$ doesn't converge.
However, there are several ways to remedy this and attach theta functions to indefinite quadratic forms.
For example one can use majorants (see \cite{siegel}) to obtain non-holomorphic modular forms.

In \cite{vigneras1,vigneras2} Vign\'eras gives a nice general construction for indefinite theta functions:
if $Q$ is a non-degenerate integer-valued quadratic form on $\Z^n$, and $p:\R^n\longrightarrow \C$ is a function that satisfies certain growth conditions and the differential equation $Dp=d\ko p$, where $d\in\Z$ and
\[D:=\ee-\frac{1}{4\pi}\Delta,\]
with $\ee:=\sum_{i=1}^n v_i \frac{\partial}{\partial v_i}$ the Euler operator and $\Delta=\Delta_Q$ the Laplace operator for $Q$, then the theta series 
\[\Theta_{Q,p} (\tau) := y^{-d/2}\sum_{\ell\in \Z^n}  p(\ell y^{1/2})\ko q^{Q(\ell)} \qquad (y=\im(\tau))\]
is a (non-holomorphic) modular form of weight $n/2+d$.
Note that both the construction for positive definite quadratic forms and Siegel's construction (\cite{siegel}) are a special case of Vign\'eras' result.

Another way to obtain indefinite theta functions is to restrict the sum over the full lattice to the sum over a cone.
In \cite{GZ} G\"ottsche and Zagier construct such indefinite theta functions for the case that the signature of $Q$ is $(n-1,1)$.
The slightly modified definition of these functions is:
let $Q$ be an integer-valued quadratic form of signature $(n-1,1)$ on $\Z^n$, let $B$ be the bilinear form associated to $Q$, let $c_1,c_2\in\overline{C}_Q:= C_Q\cup S_Q$ with $C_Q$ one of the components of $\{c\in \R^n \mid Q(c)<0\}$ and $S_Q$ the set of cusps of $C_Q$, and let $a\in R(c_1)\cap R(c_2)$, $b\in \R^n$, where $R(c)$ is $\R^n$ if $c\in C_Q$ and $\{a\in \R^n \mid B(c,a)\not\in\Z\}$ if $c\in S_Q$, then
\[\Theta_{a,b}^{c_1,c_2} (\tau):= \sum_{\ell\in a+\Z^n}\bigl\{\sign B(c_1,\ell) - \sign B(c_2,\ell)\bigr\}\ko q^{Q(\ell)}e^{2\pi i B(\ell,b)}.\]
These functions are holomorphic and their Fourier coefficients can easily be computed.
In the case that $c_1,c_2\in S_Q$ it is shown in \cite{GZ} that they are in fact modular.
For $c_1,c_2\in C_Q$ the function $\Theta_{a,b}^{c_1,c_2}$ is in general \emph{not} modular, only for special choices of $a,b$ and $c_1,c_2\in C_Q$ (see \cite{andrews} and \cite{polishchuk} for examples).\
Note that for signature $(1,1)$ such modular $\Theta_{a,b}^{c_1,c_2}$ are related to the indefinite theta functions constructed by Hecke in \cite{hecke1,hecke2}.

In \cite{zwegers} the second author showed that we can remedy the non-modularity of $\Theta_{a,b}^{c_1,c_2}$ by considering a slightly modified version:
for $c_1,c_2\in C_Q$ define
\[\widehat\Theta_{a,b}^{c_1,c_2} (\tau):= \sum_{\ell\in a+\Z^n}\Bigl\{E\Bigl(\frac{B(c_1,\ell)}{\sqrt{-Q(c_1)}}y^{1/2}\Bigr)-E\Bigl(\frac{B(c_2,\ell)}{\sqrt{-Q(c_2)}}y^{1/2}\Bigr)\Bigr\}\ko q^{Q(\ell)}e^{2\pi i B(\ell,b)},\]
where
\begin{align}\label{errorfunction}
E(z) := 2\int_0^z e^{-\pi u^2}du = \operatorname{sgn}(z) - \operatorname{sgn}(z) \int_{z^2}^\infty u^{-1/2} e^{-\pi u} du.
\end{align}
In \cite{zwegers} it is shown that $\widehat\Theta_{a,b}^{c_1,c_2}$ is modular of weight $n/2$ (alternatively, one could use the methods from \cite{vigneras1,vigneras2} to simplify the proof), but in general it is \emph{not} holomorphic.
We can view $\widehat\Theta_{a,b}^{c_1,c_2}$ as the modular ``completion'' of $\Theta_{a,b}^{c_1,c_2}$.
One application of these indefinite theta functions is that one can use them to study the modular behavior of Ramanujan's mock theta functions (see \cite{zwegers}).
Further, in certain special cases one has $\Theta_{a,b}^{c_1,c_2}=\widehat\Theta_{a,b}^{c_1,c_2}$, which explains the modularity of $\Theta_{a,b}^{c_1,c_2}$ in these cases.
Note that recently, analogous constructions have been found for quadratic forms of general signature:
see \cite{ABMP} and \cite{kudla} for signature $(n-2,2)$, and \cite{raum}, \cite{nazaroglu} and \cite{FK} (in chronological order) for the general case.

The aim of this paper is to generalize the results from \cite{zwegers} for quadratic forms of signature $(n-1,1)$ to include (spherical) polynomials.
For this we actually construct three versions for the theta function attached to a homogeneous polynomial:
a holomorphic, an almost holomorphic and a modular version.
By giving a criterion for the almost holomorphic and the modular version to agree, we obtain a construction of theta functions for quadratic forms of signature $(n-1,1)$ which are almost holomorphic modular forms.
The restriction to spherical polynomials then yields holomorphic modular forms.

\section{Definitions and statement of the main results}

For the rest of the paper we assume that the quadratic form $Q$ has signature $(n-1,1)$ and is integer-valued on $\Z^n$.
We let $A$ denote the corresponding even symmetric matrix (so $Q(v)=\frac12 v^\trans Av$) and let $B$ be the bilinear form associated to $Q$: $B(u,v)=u^\trans Av=Q(u+v)-Q(u)-Q(v)$.
Since $Q$ has signature $(n-1,1)$, the set of vectors $c\in\R^n$ with $Q(c)<0$ has two components.
If $B(c_1,c_2)<0$, then $c_1$ and $c_2$ belong to the same component, while if $B(c_1,c_2)>0$ then $c_1$ and $c_2$ belong to opposite components.
Let $C_Q$ be one of those components.
If $c_0$ is in that component, then $C_Q$ is given by:
\[C_Q :=\{ c\in \R^n \mid Q(c)<0,\ B(c,c_0)<0\}\]
We normalize the elements of $C_Q$ such that $Q(c)=-1$ and set
\[\Ce_Q:= \{ c\in \R^n \mid Q(c)=-1,\ B(c,c_0)<0\}.\]

\begin{definition}\label{defhol}
Let $f:\R^n \longrightarrow \C$ be a homogeneous polynomial of degree $d$ and let $c_1,c_2\in\Ce_Q$.
We define the holomorphic theta series associated to $Q$ and $f$ by
\[ %\Theta^{c_1,c_2}(\tau) =
\Theta^{c_1,c_2}[f](\tau) := \sum_{\ell\in\Z^n} \bigl\{\sign B(c_1,\ell) - \sign B(c_2,\ell)\bigr\}\ko f(\ell)\ko q^{Q(\ell)}.\]
\end{definition}

For the corresponding almost holomorphic and non-holomorphic versions we also need:

\begin{definition}\label{defp}
Let $\Delta=\Delta_Q:=\bigl(\frac{\partial}{\partial v}\bigr)^\trans A^{-1}\frac{\partial}{\partial v}$ denote the Laplacian associated to $Q$ (we omit $Q$ in the notation, as we take it to be fixed).
We set
\[ e^{-\Delta/8\pi}:= \sum_{k=0}^\infty \frac{(-1)^k}{(8\pi)^k k!}\ko \Delta^k,\qquad \partial_c:= c^\trans \frac{\partial}{\partial v}=\sum_{i=1}^n c_i \frac{\partial}{\partial v_i}\]
and for a homogeneous polynomial $f:\R^n\longrightarrow\C$ of degree $d$ we define $\widehat f:=e^{-\Delta /8\pi} f$ and
\[ p^c[f](v) := \sum_{k=0}^d \frac{(-1)^k}{(4\pi)^kk!}\ko E^{(k)}(B(c,v))\cdot \partial_c^k \widehat f(v).\]
\end{definition}

\begin{definition}\label{defnonhol}
Let $f:\R^n \longrightarrow \C$ be a homogeneous polynomial of degree $d$ and let $c_1,c_2\in\Ce_Q$.
We define the almost holomorphic theta series associated to $Q$ and $f$ by
\[ \widehat\Theta^{c_1,c_2}[f](\tau) := y^{-d/2} \sum_{\ell\in\Z^n} \bigl\{\sign B(c_1,\ell) - \sign B(c_2,\ell)\bigr\}\ko \widehat f(\ell y^{1/2})\ko q^{Q(\ell)}.\]
Further, we define the corresponding non-holomorphic theta series by
\[ \modtheta^{c_1,c_2}[f](\tau) := y^{-d/2} \sum_{\ell\in\Z^n} \bigl\{p^{c_1}[f](\ell y^{1/2}) - p^{c_2}[f](\ell y^{1/2})\bigr\}\ko q^{Q(\ell)}.\]
\end{definition}

\begin{remark}
(a) According to Lemma \ref{growth} $v\mapsto e^{-2\pi Q(v)} (p^{c_1}[f](v)- p^{c_2}[f](v))$ is a Schwartz function.
This ensures the absolute convergence of the sum defining $\modtheta^{c_1,c_2}[f]$.
In the proof of Lemma \ref{growth} it is further shown that for any polynomial $P$ we have
\[ \bigl| P(v)\ko e^{-2\pi Q(v)} \bigl( \sign(B(c_1,v))-\sign(B(c_2,v))\bigr)\bigr| \leq 2\ko\bigl| P(v)\bigr|\ko e^{-2\pi Q^+(v)},\]
where $Q^+$ is a positive definite quadratic from.
This directly gives us the absolute convergence of the sums in the definition of $\Theta^{c_1,c_2}[f]$ and of $\widehat\Theta^{c_1,c_2}[f]$.\\
(b) For a homogeneous polynomial $f$ of degree $d$, $\Delta^k f$ is zero (in particular for $k>d/2$) or a homogeneous polynomial of degree $d-2k$.
%(or zero for $k>d/2$).
Therefore, the degrees of the monomials in $\widehat{f}$ have the same parity as the degree of $f$, and $y^{-d/2}\widehat{f}(\ell y^{1/2})
$ is a polynomial of degree $\leq d/2$ in $1/y$.
Thus we can view $\widehat\Theta^{c_1,c_2}[f]$ as a polynomial of degree $\leq d/2$ in $1/y$ with holomorphic coefficients.
Such functions are called almost holomorphic of depth $\leq d/2$, so $\widehat\Theta^{c_1,c_2}[f]$ denotes an almost holomorphic theta series of depth $\leq d/2$.
Further, $\Theta^{c_1,c_2}[f]$ is the ``constant term'' of $\widehat\Theta^{c_1,c_2}[f]$ (viewed as a polynomial in $1/y$).
\end{remark}

We denote by $N$ the level of $A$, thus the smallest $N\in \N$ such that $NA^{-1}$ is an even matrix.
To describe the modular transformation behavior of $\modtheta^{c_1,c_2}[f]$ for a congruence subgroup of level $N$, we introduce the character $\chi$ which is defined as follows (see Theorem 2 in \cite{vigneras2}):

\begin{definition}\label{defcharacter}
Let $\bigl(\frac{\cdot}{\cdot}\bigr)$ denote the Kronecker symbol.
For $\gamma=\left(\begin{smallmatrix} a&b\\ c&d \end{smallmatrix}\right)\in\Gamma_0(N)$ we set
\[\chi(\gamma) := \Bigl(\frac{D}{d}\Bigr) \cdot v(\gamma)\quad\text{with }
\begin{cases} D=(-1)^{n/2} \det A \text{ and } v(\gamma)=1&\text{ for $n$ even,}\\
D=2\ko\det A \text{ and } v(\gamma)=\bigl(\frac{c}{d}\bigr)\, \bigl(\frac{-4}{d}\bigr)^{-n/2}&\text{ for $n$ odd.}\end{cases}\]
\end{definition}

\begin{theorem}\label{theo1}
The theta function $\modtheta^{c_1,c_2}[f]$ transforms as a modular form of weight $n/2+d$ and character $\chi$ on $\Gamma_0(N)$:
we have
\[ \modtheta^{c_1,c_2}[f]\Bigl(\frac{a\tau+b}{c\tau+d}\Bigr) = \chi(\gamma)\ko (c\tau+d)^{n/2+d}\, \modtheta^{c_1,c_2}[f](\tau)\]
for all $\gamma=\left(\begin{smallmatrix} a&b\\ c&d \end{smallmatrix}\right)\in \Gamma_0(N)$.
\end{theorem}

In order to construct almost holomorphic and holomorphic modular forms, we consider the automorphism group, which leaves the quadratic form $Q$, the lattice $\Z^n$ and the choice of the component $\Ce_Q$ unchanged:

\begin{definition}
Let
\[\aut{n}:=\bigl\{g\in\operatorname{GL}_n(\Z) \mid g^\trans A g= A,\, B(g c,c)<0\text{ for all }c\in \Ce_Q\bigr\}.\]
\end{definition}

\begin{theorem}\label{theo2}
Let $I$ be a finite set of indices.
For all $i\in I$ let $f_i$ be a homogeneous polynomial of degree $d$ and let $g_i\in\aut{n}$.
Under the assumption that
\[\sum_{i\in I} (f_i-f_i\circ g_i)=0,\]
the theta function
\[\sum_{i\in I}\modtheta^{c,g_i c}[f_i]=\sum_{i\in I}\widehat\Theta^{c,g_i c}[f_i]\]
is an almost holomorphic cusp form of weight $n/2+d$, depth $\leq d/2$ and character $\chi$ on $\Gamma_0(N)$.
Further, it doesn't depend on the choice of $c\in\Ce_Q$.
\end{theorem}

\begin{remark}
Since $\Theta^{c_1,c_2}[f]$ is the ``constant term'' of $\widehat\Theta^{c_1,c_2}[f]$, we get that the corresponding holomorphic theta function $\sum_{i\in I} \Theta^{c,g_ic}[f_i]$ is a quasimodular form (with the given weight, depth, character and subgroup).
For more details on almost holomorphic modular forms and quasimodular forms see for example section 5.3 in \cite{zagier}.
\end{remark}

\begin{definition}
We call a polynomial $f:\R^n\longrightarrow\C$ \emph{spherical} (of degree $d$) if it is homogeneous (of degree $d$) and vanishes under the Laplacian, i.\,e. $\Delta f=0$.
\end{definition}

\begin{remark}\label{remsph}
If $f$ is spherical of degree $d$, then we have $\widehat f=e^{-\Delta /8\pi} f=f$ and $y^{-d/2}\widehat f(\ell y^{1/2})= y^{-d/2} f(\ell y^{1/2})=f(\ell)$.
Hence the holomorphic theta function $\Theta^{c_1,c_2}[f]$ and the almost holomorphic theta function $\widehat\Theta^{c_1,c_2}[f]$ agree.
This observation immediately leads to the following corollary to Theorem \ref{theo2}.
\end{remark}

\begin{corollary}\label{cor1}
Let $I$, $f_i$ and $g_i$ be as in Theorem \ref{theo2}, with the additional condition that $f_i$ is spherical for all $i\in I$.
Under the assumption that
\[\sum_{i\in I} (f_i-f_i\circ g_i)=0,\]
the theta function
\[\sum_{i\in I}\modtheta^{c,g_i c}[f_i]=\sum_{i\in I}\Theta^{c,g_i c}[f_i]\]
is a (holomorphic) cusp form of weight $n/2+d$ and character $\chi$ on $\Gamma_0(N)$.
Further, it doesn't depend on the choice of $c\in\Ce_Q$.
\end{corollary}

\begin{remark}\label{remexamples}
Since $E$ is odd, we get $p^c[f](-v)=(-1)^{d+1}p^c[f](v)$.
Hence if the degree $d$ of $f$ is even we trivially have
that $\modtheta^{c_1,c_2}[f]$ is identically 0.
Similarly $\widehat\Theta^{c_1,c_2}[f]$ and $\Theta^{c_1,c_2}[f]$ also vanish.
In this case, non-trivial results can still be obtained by introducing characteristics $a,b\in\Q^n$ and setting
\[ \modtheta_{a,b}^{c_1,c_2}[f](\tau) := y^{-d/2} \sum_{\ell\in a+\Z^n} \bigl\{p^{c_1}[f](\ell y^{1/2}) - p^{c_2}[f](\ell y^{1/2})\bigr\}\ko q^{Q(\ell)}\ko e^{2\pi iB(\ell,b)},\]
and similarly for the holomorphic and the almost holomorphic versions (as was done in \cite{zwegers}, where the case $d=0$ is considered).
Analogously, one can include periodic functions on $\Z^n$, that is, functions $m:\Z^n\longrightarrow \C$ such that there is an $L\in\N$ for which we have $m(\ell +\ell')=m(\ell)$ for all $\ell\in\Z^n$ and all $\ell'\in L\Z^n$, and consider
\[\begin{split}
\Theta^{c_1,c_2}[m,f](\tau) &:= \sum_{\ell\in\Z^n} \bigl\{\sign B(c_1,\ell) - \sign B(c_2,\ell)\bigr\}\ko m(\ell)\ko f(\ell)\ko q^{Q(\ell)},\\
\widehat\Theta^{c_1,c_2}[m,f](\tau) &:= y^{-d/2}
 \sum_{\ell\in\Z^n} \bigl\{\sign B(c_1,\ell) - \sign B(c_2,\ell)\bigr\}\ko m(\ell)\ko \widehat f(\ell y^{1/2})\ko q^{Q(\ell)},\\
\modtheta^{c_1,c_2}[m,f](\tau) &:= y^{-d/2} \sum_{\ell\in\Z^n} \bigl\{p^{c_1}[f](\ell y^{1/2}) - p^{c_2}[f](\ell y^{1/2})\bigr\}\ko m(\ell)\ko q^{Q(\ell)}.
\end{split}\]
With a slight generalization of Vign\'eras' result one can show that $\modtheta^{c_1,c_2}[m,f]$ again transforms as a modular form of weight $n/2+d$, but the subgroup and the character now also depend on the choice of $m$.
We omit the details.
Theorem \ref{theo2} and Corollary \ref{cor1} then generalize to:
for all $i\in I$ let $m_i$ be a periodic function on $\Z^n$, let $f_i$ be a homogeneous polynomial of degree $d$, let $\widetilde f_i=m_i\cdot f_i$ and let $g_i\in\aut{n}$.
Under the assumption that $\sum_{i\in I} (\widetilde f_i-\widetilde f_i\circ g_i)=0$, the theta function
\[\sum_{i\in I}\modtheta^{c,g_i c}[m_i,f_i]=\sum_{i\in I}\widehat\Theta^{c,g_i c}[m_i,f_i]\]
is an almost holomorphic cusp form of weight $n/2+d$ and depth $\leq d/2$.
If we further assume that $f_i$ is spherical for all $i\in I$, then
\[\sum_{i\in I}\modtheta^{c,g_i c}[m_i,f_i]=\sum_{i\in I}\Theta^{c,g_i c}[m_i,f_i]=\sum_{i\in I}\Theta^{c,g_i c}[\widetilde f_i]\]
is a (holomorphic) cusp form of weight $n/2+d$.
%The same approach can be applied if the automorphism group $\aut$ contains few elements of infinite order, this will be carried out explicitly in some of the following examples.
\end{remark}

\section{Proof of Theorems \ref{theo1} and \ref{theo2}, and Corollary \ref{cor1}}

To prove Theorem \ref{theo1} we'll use the results from \cite{vigneras1,vigneras2}, for which we have to show that the function $p:=p^{c_1}[f]-p^{c_2}[f]$ satisfies the necessary growth conditions (Lemma \ref{growth}) and the differential equation $Dp=d\ko p$ (Lemma \ref{lemdif}).

\begin{lemma}\label{growth}
For $c_1,c_2\in\Ce_Q$ and $f$, $p^c[f]$ as in Definition \ref{defp}
\[v\mapsto e^{-2\pi Q(v)} \bigl(p^{c_1}[f](v)- p^{c_2}[f](v)\bigr)\]
is a Schwartz function.
\end{lemma}

\begin{proof}
If $c_1$ and $c_2$ are linearly dependent, the expression $p^{c_1}[f](v)- p^{c_2}[f](v)$ vanishes, thus we assume that they are linearly independent.

Since $p^c[f]$ is the finite sum of $C^\infty$-functions, $v\mapsto e^{-2\pi Q(v)} (p^{c_1}[f](v)- p^{c_2}[f](v))$ is also a $C^\infty$-function.
To show that it is a Schwartz function we begin by splitting $p^c[f](v)$ as the sum of $E(B(c,v))\ko \widehat f(v)$ and
\begin{equation}\label{defptilde}
\widetilde p^c[f](v):=\sum_{k=1}^d \frac{(-1)^k}{(4\pi)^kk!}\ko E^{(k)}(B(c,v))\cdot \partial_c^k \widehat f(v).
\end{equation}
We have $E'(B(c,v))=2e^{-\pi B(c,v)^2}$ and with induction we can easily see that for all $k\in\N$ we can write $E^{(k)}(B(c,v))$ as a polynomial (in $v$) times $e^{-\pi B(c,v)^2}$.
Hence we can write $e^{-2\pi Q(v)}\widetilde p^c[f](v)$ as a polynomial times $e^{-2\pi Q_c(v)}$, where $Q_c(v):= Q(v)+\frac12B(c,v)^2$.
Since $Q_c$ is a positive definite quadratic form (Lemma 2.5 in \cite{zwegers}), $v\mapsto e^{-2\pi Q_c(v)}$ is a Schwartz function and hence so is $v\mapsto e^{-2\pi Q(v)}\widetilde p^c[f](v)$.
What remains to be shown is that $v\mapsto e^{-2\pi Q(v)}[E(B(c_1,v))-E(B(c_2,v))]\ko \widehat f(v)$ is a Schwartz function.
Since $\widehat f$ is a polynomial we actually only need to show that $v\mapsto e^{-2\pi Q(v)} [E(B(c_1,v))-E(B(c_2,v))]$ is one.
By induction on the total number of derivatives we can easily see that any higher-order partial derivative of $e^{-2\pi Q(v)}[E(B(c_1,v))-E(B(c_2,v))]$ is of the form
\[ P_0(v)\ko e^{-2\pi Q(v)} \bigl( E(B(c_1,v))-E(B(c_2,v))\bigr) + P_1(v)\ko e^{-2\pi Q(v)} e^{-\pi B(c_1,v)^2} - P_2(v)\ko e^{-2\pi Q(v)} e^{-\pi B(c_2,v)^2},\]
where $P_0$, $P_1$ and $P_2$ are polynomials.
By the same argument as before $P_i(v)\ko e^{-2\pi Q(v)} e^{-\pi B(c_i,v)^2}$ is a Schwartz function, so it suffices to show that for any polynomial $P$
\begin{align}\label{est}
\bigl| P(v)\ko e^{-2\pi Q(v)} \bigl( E(B(c_1,v))-E(B(c_2,v))\bigr)\bigr|
\end{align}
is bounded on $\R^n$.
For this we use more or less the same arguments as in \cite{zwegers}:
rewriting $E$ as in (\ref{errorfunction}), an upper bound for the expression in \eqref{est} is given by the sum of the three expressions
\begin{equation}\label{est_int}
\bigl| P(v)\ko e^{-2\pi Q(v)} \sign(B(c_i,v))\ko \beta(B(c_i,v)^2) \bigr|, \qquad (i=1,2)
\end{equation}
where $\beta(x):=\int_x^\infty u^{-1/2} e^{-\pi u} du$, and
\begin{equation}\label{est_diff}
\bigl| P(v)\ko e^{-2\pi Q(v)} \bigl( \sign(B(c_1,v))-\sign(B(c_2,v))\bigr)\bigr|.
\end{equation}
Considering \eqref{est_int}, we use $0\leq \beta(x) \leq e^{-\pi x}$ for $x\in \R_{\geq 0}$ to find
\begin{equation}\label{estintq}
\bigl| P(v)\ko e^{-2\pi Q(v)} \sign(B(c_i,v))\ko \beta(B(c_i,v)^2) \bigr| \leq \bigl| P(v)\bigr|\ko e^{-2\pi Q_{c_i}(v)},
\end{equation}
where the last expression is bounded on $\R^n$ (as before).

Obviously, \eqref{est_diff} vanishes if $\sign(B(c_1,v))\ko\sign(B(c_2,v))>0$ holds.
For linearly independent vectors $c_1,c_2\in \Ce_Q$ we can check (Lemma 2.6 in \cite{zwegers}) that $Q^+(v):=Q(v)+\frac{B(c_1,c_2)}{4-B(c_1,c_2)^2}B(c_1,v)\ko B(c_2,v)$ is a positive definite quadratic form.
As $c_1,c_2\in \Ce_Q$, we have $B(c_1,c_2)<0$ and $c_1,c_2$ span a subspace where $Q$ is of signature $(1,1)$, thus the determinant of $\left(\begin{smallmatrix} 2Q(c_1)&B(c_1,c_2)\\ B(c_1,c_2)&2Q(c_2) \end{smallmatrix}\right)$ is negative, i.\,e. $4-B(c_1,c_2)^2<0$.
Hence, for $\sign(B(c_1,v))\ko\sign(B(c_2,v))\leq 0$ we have $Q(v)\geq Q^+(v)$ and so
\[\bigl| P(v)\ko e^{-2\pi Q(v)} \bigl( \sign(B(c_1,v))-\sign(B(c_2,v))\bigr)\bigr|\leq 2\ko\bigl| P(v)\bigr|\ko e^{-2\pi Q^+(v)},\]
where $2\ko | P(v)|\ko e^{-2\pi Q^+(v)}$ is bounded on $\R^n$.

Thus we have shown that the expression in (\ref{est}) is bounded, which completes the proof.
\end{proof}

\begin{lemma}\label{lemdif}
For $f$, $\widehat f$ and $p^c[f]$ as in Definition \ref{defp} we have $D\widehat f=d\ko \widehat f$ and $Dp^c[f]=d\ko p^c[f]$.
\end{lemma}

For the proof we need:

\begin{lemma}\label{commu}
We have $D\ko e^{-\Delta/8\pi}= e^{-\Delta/8\pi}\ee$.
Further, the differential operators $D$ and $\partial_c^k$ ($k\in\N_0$) satisfy the commutator relation $[D,\partial_c^k]=-k\partial_c^k$.
\end{lemma}

\begin{proof}
One can easily check that $[\ee,\Delta]=-2 \Delta$ and $[D,\partial_c]=-\partial_c$ hold, which by induction generalize directly to $[\ee,\Delta^k]=-2k\Delta^k$ and $[D,\partial_c^k]=-k\partial_c^k$ for all $k\in\N_0$.
Further,
\[ [\ee, e^{-\Delta/8\pi}]= \sum_{k=0}^\infty \frac{(-1)^k}{(8\pi)^k k!}\ko [\ee,\Delta^k]= \sum_{k=0}^\infty \frac{(-1)^k}{(8\pi)^k k!}\ko(-2k\Delta^k )= \frac1{4\pi} \sum_{k=0}^\infty \frac{(-1)^k}{(8\pi)^k k!}\ko \Delta^{k+1}=\frac\Delta{4\pi}e^{-\Delta/8\pi}\]
and so $D\ko e^{-\Delta/8\pi}= e^{-\Delta/8\pi}\ee$.
\end{proof}

\begin{proofof}{Lemma \ref{lemdif}}
Since $f$ is homogeneous of degree $d$ it satisfies $\ee f=d\ko f$ and so Lemma \ref{commu} gives
\[ D\widehat f = D\ko e^{-\Delta/8\pi}f= e^{-\Delta/8\pi}\ee f= d\ko e^{-\Delta/8\pi}f = d\ko \widehat f.\]
Further, we have
\[ D\ko \partial_c^k \widehat f= \partial_c^k D \widehat f - k\ko \partial_c^k \widehat f= (d-k)\ko \partial_c^k\widehat f\]
and a direct computation gives
\[ \ee \bigl(E(B(c,v))\bigr)= 2 B(c,v)\ko e^{-\pi B(c,v)^2} \qquad \text{and} \qquad \Delta \bigl(E(B(c,v))\bigr) = 8\pi B(c,v)\ko e^{-\pi B(c,v)^2}.\]
Hence we find $D(E(B(c,v)))=0$ and
\[ D\ko \partial_c^k\bigl(E(B(c,v))\bigr)= \partial_c^kD\bigl(E(B(c,v))\bigr)-k\ko \partial_c^k\bigl(E(B(c,v))\bigr)=-k\ko \partial_c^k\bigl(E(B(c,v))\bigr).\]
Since $\partial_c^k\bigl(E(B(c,v))\bigr)=(-2)^k E^{(k)} (B(c,v))$, this yields
\[ D\bigl(E^{(k)} (B(c,v))\bigr) = -k E^{(k)} (B(c,v)).\]
For the product of two functions we have
\[ D(f_1 \cdot f_2) = Df_1 \cdot f_2 +f_1\cdot Df_2 -\frac1{2\pi} \sum_{i=1}^n\sum_{j=1}^n (A^{-1})_{ij} \frac{\partial f_1}{\partial v_i} \frac{\partial f_2}{\partial v_j}.\]
Setting $f_1(v)=E^{(k)}(B(c,v))$ and $f_2=\partial_c^k\widehat f $ gives
\[\begin{split}
\sum_{i=1}^n\sum_{j=1}^n (A^{-1})_{ij}\ko \frac{\partial f_1}{\partial v_i} \frac{\partial f_2}{\partial v_j} &= \sum_{i=1}^n\sum_{j=1}^n (A^{-1})_{ij} (Ac)_i\ko E^{(k+1)} (B(c,v))\ko \frac{\partial}{\partial v_j}\ko \partial_c^k \widehat f(v)\\
&= \sum_{j=1}^n E^{(k+1)} (B(c,v))\, c_j\frac{\partial}{\partial v_j}\ko \partial_c^k \widehat f(v)= E^{(k+1)}(B(c,v))\cdot \partial_c^{k+1} \widehat f(v)
\end{split}\]
and so
\[D\bigl(E^{(k)}(B(c,v))\cdot \partial_c^k \widehat f(v)\bigr)= (d-2k)\ko E^{(k)} (B(c,v))\cdot \partial_c^k \widehat f(v) -\frac1{2\pi}E^{(k+1)}(B(c,v))\cdot \partial_c^{k+1} \widehat f(v).\]
Hence
\[ \begin{split}
D\ko p^c[f](v)&= \sum_{k=0}^d \frac{(-1)^k}{(4\pi)^k k!}\ko \bigl\{(d-2k) \ko E^{(k)} (B(c,v))\cdot \partial_c^k \widehat f(v) -\frac1{2\pi}E^{(k+1)}(B(c,v))\cdot \partial_c^{k+1} \widehat f(v) \bigr\}\\
&= \sum_{k=0}^d (d-2k)\frac{(-1)^k}{(4\pi)^k k!}\ko E^{(k)} (B(c,v))\cdot \partial_c^k \widehat f(v) +\frac1{2\pi}\sum_{k=1}^d \frac {(-1)^k}{(4\pi)^{k-1} (k-1)!}\ko E^{(k)} (B(c,v))\cdot \partial_c^k \widehat f(v)\\
&= d \sum_{k=0}^d \frac{(-1)^k}{(4\pi)^k k!}\ko E^{(k)} (B(c,v))\cdot \partial_c^k \widehat f(v) = d\ko p^c[f](v),
\end{split}\]
where we have used that $\partial_c^{d+1} \widehat f(v)=0$.
To prove this last identity we observe that if $f$ is homogeneous of degree $d$, then $\partial_c f$ is homogeneous of degree $d-1$.
Hence $\partial_c^{d+1}f=0$ and so
\[\partial_c^{d+1}\widehat f = \partial_c^{d+1} e^{-\Delta/8\pi}f= e^{-\Delta/8\pi} \partial_c^{d+1}f=0.\]
This finishes the proof.
\end{proofof}

\begin{proofof}{Theorem \ref{theo1}}
We set $p:=p^{c_1}[f]- p^{c_2}[f]$.
According to Lemma \ref{growth} $v\mapsto e^{-2\pi Q(v)} p(v)$ is a Schwartz function.
Further, it follows directly from Lemma \ref{lemdif} that $p$ satisfies the differential equation $Dp=d\ko p$.
Hence we can apply Theorems 1 and 2 from \cite{vigneras2} to get the desired result.
\end{proofof}

As is usual in the theory of theta functions, we additionally introduce the characteristic $\lambda$ in the dual lattice $(\Z^n)^*=A^{-1} \Z^n$ to be able to study the modular transformation properties of $\modtheta^{c_1,c_2}[f]$ in more detail.

\begin{definition}
Let $\lambda\in A^{-1}\Z^n$ and let $f$ and $c_1,c_2$ be as in Definitions \ref{defhol} and \ref{defnonhol}.
We define
\[\begin{split}
\Theta_\lambda^{c_1,c_2}[f](\tau) &:= \sum_{\ell\in\lambda+\Z^n} \bigl\{\sign B(c_1,\ell) - \sign B(c_2,\ell)\bigr\}\ko f(\ell)\ko q^{Q(\ell)},\\
\widehat\Theta_\lambda^{c_1,c_2}[f](\tau) &:= y^{-d/2} \sum_{\ell\in\lambda+\Z^n} \bigl\{\sign B(c_1,\ell) - \sign B(c_2,\ell)\bigr\}\ko \widehat f(\ell y^{1/2})\ko q^{Q(\ell)},\\
\modtheta_\lambda^{c_1,c_2}[f](\tau) &:= y^{-d/2} \sum_{\ell\in\lambda+\Z^n} \bigl\{p^{c_1}[f](\ell y^{1/2}) - p^{c_2}[f](\ell y^{1/2})\bigr\}\ko q^{Q(\ell)}.
\end{split}\]
\end{definition}

\begin{remark}
Since these definitions depend only on $\lambda$ modulo $\Z^n$, we will consider $\lambda$ to be in the finite set $A^{-1}\Z^n/\Z^n$.
\end{remark}

The modular transformation properties for the generators $T$ and $S$ of $\Gamma_1:=\operatorname{SL}_2(\Z)$ we immediately get from (1) and (3) in \cite{vigneras2}.
Note that there is a typo in (3), which we have corrected here ($e^{2\pi iB(\lambda,\mu)}$ is missing).

\begin{lemma}[Vign\'eras \cite{vigneras2}]\label{lemtransTS}
The theta functions with characteristic satisfy
\[ \begin{split}
\modtheta_\lambda^{c_1,c_2}[f](\tau+1) &= e^{2\pi iQ(\lambda)}\ko \modtheta_\lambda^{c_1,c_2}[f](\tau),\\
\modtheta_\lambda^{c_1,c_2}[f](-1/\tau) &= (-i\tau)^{n/2+d}\ko \frac{(-i)^{d+1}}{\sqrt{|\det A|}} \sum_{\mu \in A^{-1}\Z^n/\Z^n} e^{2\pi iB(\lambda,\mu)}\ko \modtheta_\mu^{c_1,c_2}[f](\tau).
\end{split}\]
\end{lemma}

Using these we'll show:

\begin{lemma}\label{lemmod}
For all $\gamma\in\Gamma_1$ we can write $\modtheta^{c_1,c_2}[f]\big|_{n/2+d}\gamma$ as
\[\sum_{\lambda \in A^{-1}\Z^n/\Z^n} \varphi_\gamma(\lambda)\ko \modtheta_\lambda^{c_1,c_2}[f],\]
where $\varphi_\gamma: A^{-1}\Z^n/\Z^n \longrightarrow \C$ satisfies $\varphi_\gamma \circ g = \varphi_\gamma$ for all $g\in\aut{n}$.
\end{lemma}

\begin{remark}\label{rembij}
(a) We can easily check that left-multiplication with $g\in\aut{n}$ is a bijection from $A^{-1}\Z^n/\Z^n$ into itself.
Hence for a function $\varphi$ on $A^{-1}\Z^n/\Z^n$, $\varphi\circ g$ is also a function on $A^{-1}\Z^n/\Z^n$.\\
(b) In \cite{schoeneberg} and \cite{shimura} explicit formulas are given for the coefficients $\varphi_\gamma(\lambda)$ for the case that $Q$ is positive definite.
Similar formulas hold for the indefinite case, but these are not given explicitly in \cite{vigneras2}.
However, for our purposes it suffices to know that the coefficients $\varphi_\gamma$ satisfy $\varphi_\gamma \circ g = \varphi_\gamma$ for all $g\in\aut{n}$.
\end{remark}

\begin{proofof}{Lemma \ref{lemmod}}
We say that a function is of the right form if we can write it as
\[\sum_{\lambda \in A^{-1}\Z^n/\Z^n} \varphi(\lambda)\ko \modtheta_\lambda^{c_1,c_2}[f],\]
where $\varphi \circ g = \varphi$ holds for all $g\in\aut{n}$.
We first observe that $\modtheta^{c_1,c_2}[f]$ is of the right form:
we have
\[ \modtheta^{c_1,c_2}[f] = \sum_{\lambda \in A^{-1}\Z^n/\Z^n} \varphi(\lambda)\ko \widehat\Theta_\lambda^{c_1,c_2}[f],\]
with $\varphi(\lambda) =1$ if $\lambda \equiv 0 \smod{\Z^n}$ and 0 otherwise.
Indeed we have $\varphi\circ g=\varphi$ for all $g\in\aut{n}$.
Further, it follows from Lemma \ref{lemtransTS} that if $h$ is of the right form, then so are $h\big|_{n/2+d} T$ and $h\big|_{n/2+d} S$:
\[\begin{split}
h\big|_{n/2+d} T &= \sum_{\lambda \in A^{-1}\Z^n/\Z^n} \varphi(\lambda)\ko \modtheta_\lambda^{c_1,c_2}[f]\big|_{n/2+d} T \\
&= \sum_{\lambda \in A^{-1}\Z^n/\Z^n} \varphi(\lambda)\ko e^{2\pi iQ(\lambda)}\ko \modtheta_\lambda^{c_1,c_2}[f]=\sum_{\lambda \in A^{-1}\Z^n/\Z^n} \varphi_1(\lambda)\ko \modtheta_\lambda^{c_1,c_2}[f],\\
h\big|_{n/2+d} S &=\sum_{\lambda \in A^{-1}\Z^n/\Z^n} \varphi(\lambda)\ko \modtheta_\lambda^{c_1,c_2}[f]\big|_{n/2+d} S \\
&= \sum_{\lambda \in A^{-1}\Z^n/\Z^n} \varphi(\lambda) \frac{(-i)^{n/2+2d+1}}{\sqrt{|\det A|}} \sum_{\mu \in A^{-1}\Z^n/\Z^n} e^{2\pi iB(\lambda,\mu)}\ko \modtheta_\mu^{c_1,c_2}[f]= \sum_{\mu \in A^{-1}\Z^n/\Z^n} \varphi_2(\mu)\ko \modtheta_\mu^{c_1,c_2}[f],
\end{split}\]
with $\varphi_1(\lambda):= \varphi(\lambda)\ko e^{2\pi iQ(\lambda)}$ and
\[\varphi_2(\mu) :=  \frac{(-i)^{n/2+2d+1}}{\sqrt{|\det A|}} \sum_{\lambda \in A^{-1}\Z^n/\Z^n}e^{2\pi iB(\lambda,\mu)}\ko \varphi(\lambda).\]
Since we assume $\varphi\circ g =\varphi$ for all $g\in\aut{n}$, it follows directly that we also have $\varphi_i\circ g=\varphi_i$ ($i=1,2$) for all $g\in\aut{n}$ (for $i=2$ replace $(\lambda,\mu)$ by $(g\lambda,g\mu)$ and use part (a) of Remark \ref{rembij}).

Since  the group $\Gamma_1$ is generated by $T$ and $S$, and $\modtheta^{c_1,c_2}[f]$ is of the right form, it now follows that $\modtheta^{c_1,c_2}[f]\big|_{n/2+d}\gamma$ is of the right form for all $\gamma\in \Gamma_1$.
\end{proofof}

\begin{remark}\label{remdep}
From these computations we also see that the coefficients $\varphi_\gamma(\lambda)$ don't depend on the choice of $c_1$, $c_2$.
Further, they do depend on the degree $d$ of $f$, but not on $f$ itself.
\end{remark}

\begin{proofof}{Theorem \ref{theo2}}
The key to the proof is that under the assumption that $\sum_{i\in I}(f_i-f_i\circ g_i)=0$ is satisfied, the modular theta series $\sum_{i\in I} \modtheta^{c,g_i c}[f_i]$ and the almost holomorphic theta series $\sum_{i\in I}\widehat\Theta^{c,g_i c}[f_i]$ coincide.
In fact, we will show that more generally we have
\begin{equation}\label{mainid}
\sum_{i\in I,\,\lambda \in A^{-1}\Z^n/\Z^n} \varphi(\lambda)\ko \modtheta_\lambda^{c,g_i c}[f_i]=\sum_{i\in I,\, \lambda \in A^{-1}\Z^n/\Z^n} \varphi(\lambda)\ko \widehat\Theta_\lambda^{c,g_i c}[f_i]
\end{equation}
if $\varphi \circ g = \varphi$ holds for all $g\in\aut{n}$.
As in the proof of Lemma \ref{growth}, we split $p^c[f](v)$ as the sum of $E(B(c,v))\ko \widehat f(v)$ and $\widetilde p^c[f](v)$.
Further, we split $E$ by using \eqref{errorfunction}.
Together this gives
\begin{equation}\label{splitp}
p^{c_1}[f](v) - p^{c_2}[f](v)=\bigl\{\sign B(c_1,v) - \sign B(c_2,v)\bigr\} \widehat f(v) +\breve{p}^{c_1}[f](v)-\breve{p}^{c_2}[f](v),
\end{equation}
where
\[ \breve{p}^c[f](v)=\widetilde p^{c}[f](v)-\sign(B(c,v))\ko \beta(B(c,v)^2)\ko \widehat f(v).\]
In the proof of Lemma \ref{growth} we have seen that $v\mapsto e^{-2\pi Q(v)}\widetilde p^c[f](v)$ is a Schwartz function and from \eqref{estintq} we get
\[ \bigl| e^{-2\pi Q(v)}\ko \sign(B(c,v))\ko \beta(B(c,v)^2)\ko \widehat f(v) \bigr| \leq \bigl| \widehat f(v)\bigr|\ko e^{-2\pi Q_c(v)},\]
where $Q_c$ is a positive definite quadratic form.
Hence for $c\in\Ce_Q$
\[ \theta_\lambda^c[f](\tau):= y^{-d/2} \sum_{\ell\in\lambda+\Z^n} \breve p^c[f](\ell y^{1/2}) \ko q^{Q(\ell)}\]
converges absolutely.
With \eqref{splitp} we thus obtain
\begin{equation}\label{splittheta}
\modtheta_\lambda^{c_1,c_2}[f] = \widehat\Theta_\lambda^{c_1,c_2}[f] + \theta_\lambda^{c_1}[f] - \theta_\lambda^{c_2}[f].
\end{equation}
Now let $g\in \aut{n}$.
We can easily check $\Delta(f\circ g)= (\Delta f)\circ g$, which gives $\widehat f\circ g= \widehat{f\circ g}$, and $(\partial_{gc}^k \widehat f)\circ g=\partial_c^k(\widehat f\circ g)$.
Using these in the definition of $\widetilde p^c[f]$ (equation \eqref{defptilde}) and of $\breve p^c[f]$ we find $\widetilde p^{gc}[f](gv)=\widetilde p^c[f\circ g](v)$ and $\breve p^{gc}[f](gv)=\breve p^c[f\circ g](v)$.
Replacing $(c,\lambda,\ell)$ by $(gc,g\lambda,g\ell)$ in the definition of $\theta_\lambda^c[f]$ then gives
\begin{equation}\label{phiaut}
\theta_{g\lambda}^{gc}[f]= \theta_\lambda^c[f\circ g].
\end{equation}
Using \eqref{splittheta} we find
\[ \begin{split}
\sum_{i\in I,\,\lambda \in A^{-1}\Z^n/\Z^n} \varphi(\lambda)\ko \modtheta_\lambda^{c,g_i c}[f_i]&=\sum_{i\in I,\, \lambda \in A^{-1}\Z^n/\Z^n} \varphi(\lambda)\ko \widehat\Theta_\lambda^{c,g_i c}[f_i] \\
&\qquad + \sum_{i\in I,\, \lambda \in A^{-1}\Z^n/\Z^n} \varphi(\lambda)\ko \theta_\lambda^c[f_i] -  \sum_{i\in I,\, \lambda \in A^{-1}\Z^n/\Z^n} \varphi(\lambda)\ko \theta_\lambda^{g_ic}[f_i].
\end{split}\]
Replacing $\lambda$ by $g_i\lambda$, assuming $\varphi \circ g = \varphi$ for all $g\in\aut{n}$, and using \eqref{phiaut} we get
\[ \sum_{\lambda \in A^{-1}\Z^n/\Z^n} \varphi(\lambda)\ko \theta_\lambda^{g_ic}[f_i]=\sum_{\lambda \in A^{-1}\Z^n/\Z^n} \varphi(g_i\lambda)\ko \theta_{g_i\lambda}^{g_ic}[f_i]=\sum_{\lambda \in A^{-1}\Z^n/\Z^n} \varphi(\lambda)\ko \theta_\lambda^{c}[f_i\circ g_i]\]
and so we obtain
\[ \sum_{i\in I,\, \lambda \in A^{-1}\Z^n/\Z^n} \varphi(\lambda)\ko \theta_\lambda^c[f_i] -  \sum_{i\in I,\, \lambda \in A^{-1}\Z^n/\Z^n} \varphi(\lambda)\ko \theta_\lambda^{g_ic}[f_i]=\sum_{\lambda \in A^{-1}\Z^n/\Z^n} \varphi(\lambda)\ko \theta_\lambda^c\Bigl[\sum_{i\in I} (f_i-f_i\circ g_i)\Bigr]=0,\]
which proves \eqref{mainid}.

We have already seen in the proof of Lemma \ref{lemmod} that we can choose $\varphi$ as $\varphi(\lambda) =1$ if $\lambda \equiv 0 \smod{\Z^n}$ and 0 otherwise, which gives that the modular theta series $\sum_{i\in I} \modtheta^{c,g_i c}[f_i]$ and the almost holomorphic theta series $\sum_{i\in I}\widehat\Theta^{c,g_i c}[f_i]$ agree, thus we obtain an almost holomorphic theta series of depth $\leq d/2$ with the desired transformation behavior.

In the Fourier expansion of $\widehat\Theta_\lambda^{c_1,c_2}[f]$ only positive powers of $q$ occur as the quadratic form $Q$ is bounded from below by a positive definite quadratic form $Q^{+}$ on the support of
$\bigl\{\sign B(c_1,\ell) - \sign B(c_2,\ell)\bigr\}\ko \widehat f(\ell y^{1/2})$ (as shown in the proof of Lemma \ref{growth}).
Using \eqref{mainid}, Lemma \ref{lemmod}, Remark \ref{remdep} and again \eqref{mainid} we can write for any $\gamma\in \Gamma_1$ the function $\sum_{i\in I}\widehat\Theta^{c,g_i c}[f_i]\big|_{n/2+d}\gamma$ as a linear combination of $\widehat\Theta^{c,g_ic}_\lambda[f_i]$ and thus we have a Fourier expansion with positive powers of $q$ in any cusp (where the Fourier coefficients are polynomials of degree $\leq d/2$ in $1/y$).
This shows that $\sum_{i\in I} \modtheta^{c,g_i c}[f_i]=\sum_{i\in I}\widehat\Theta^{c,g_i c}[f_i]$ is an almost holomorphic cusp form of depth $d/2$ with the given modular transformation properties.

Exactly as in the proof of $\widetilde p^{gc}[f](gv)=\widetilde p^c[f\circ g](v)$ we also have $p^{gc}[f](gv)= p^c[f\circ g](v)$.
Replacing $(c_1,c_2,\ell)$ by $(gc_1,gc_2,g\ell)$ in the definition of $\modtheta^{c_1,c_2}[f]$ hence gives
\[ \modtheta^{gc_1,gc_2} [f] = \modtheta^{c_1,c_2} [f\circ g]\]
for all $g\in\aut{n}$, and so for $c,c'\in\Ce_Q$ we have
\[ \modtheta^{c,g_i c}[f_i] -  \modtheta^{c'\hspace{-.26em},g_i c'}[f_i] = \modtheta^{c,c'}[f_i] -  \modtheta^{g_ic,g_i c'}[f_i] =   \modtheta^{c,c'}[f_i] -  \modtheta^{c,c'}[f_i\circ g_i] = \modtheta^{c,c'}[f_i-f_i\circ g_i],\]
where in the first step we used the trivial identity
\[ \bigl(p^c[f_i]-p^{g_ic}[f_i]\bigr) - \bigl(p^{c'}[f_i]-p^{g_ic'}[f_i]\bigr)=\bigl(p^c[f_i]-p^{c'}[f_i]\bigr) - \bigl(p^{g_i c}[f_i]-p^{g_ic'}[f_i]\bigr).\]
Summing over all $i\in I$ then gives
\[ \sum_{i\in I} \modtheta^{c,g_i c}[f_i]=\sum_{i\in I} \modtheta^{c'\hspace{-.26em},g_i c'}[f_i],\]
which proves the last part of the theorem.
\end{proofof}

\begin{proofof}{Corollary \ref{cor1}}
The result follows immediately from Theorem \ref{theo2} and Remark \ref{remsph}.
\end{proofof}

\section{Explicit examples}
Finally, we give explicit examples which can be constructed by using the main results of this work.
We already obtain a vast variety of nice examples for low dimensions.
Considering some quadratic forms more thoroughly, we can also make statements about the number of different modular forms we might get (see Example \ref{ex284}) and give a general construction for specific quadratic forms of level $4N$ (see Examples \ref{ex2N} and \ref{exN}).
If possible, we identify the theta series as eta quotients and define as usual $\eta(\tau):=q^{1/24}\prod_{n=1}^\infty (1-q^n)$ and $\eta_M(\tau):=\eta(M\tau)$.
Note that we don't actually prove these identities, we only verified them by computing the first 1000 coefficients in the Fourier series.
For the proof one would have to use results on eta quotients (as for example in \cite{GH}) to determine the exact modular transformation behavior% of the modular forms with regard to the congruence subgroups of $\Gamma_1$ on which the modular forms transform like modular forms as well as the character, but we will not do this here
, as well as determine the corresponding Sturm bound (see \cite{sturm} for modular forms of integral weight and \cite{KP} for modular forms of half-integral weight).

Studying binary quadratic forms of signature $(1,1)$, we do not seem to obtain any interesting examples when we include spherical polynomials of degree $d>0$.
%which do not appear as derivatives of the constructions for constants. 
For $d=0$ though, we can introduce for $i\in I$ periodic functions $m_i$ and choose $g_i\in \aut{2}$ and polynomials $f_i$ such that the functions $\widetilde{f}_i=m_i\cdot f_i$ fulfill the assumption $\sum_{i\in I}(\widetilde{f}_i-\widetilde{f}_i\circ g_i)=0 $ in Remark \ref{remexamples}.
We begin this section by giving some very simple examples for this case choosing $I=\{1\}$ and $f_1\equiv 1$.

\begin{example}
Let $Q(v)=v_1^2+5v_1v_2+v_2^2$.
We take $g=\left(\begin{smallmatrix}5&1\\ -1&0 \end{smallmatrix}\right)\in \aut{2}$ and $c=\frac{1}{\sqrt{21}}\left(\begin{smallmatrix}-2\\ 5 \end{smallmatrix}\right)\in \Ce_Q$.
As a periodic function $m:\Z^2\longrightarrow\C$ with $m\circ g=m$ we choose \[m(v)=\left(\frac{-3}{v_1+v_2}\right),\] where $\left(\frac{-3}{\cdot}\right)$ is odd and has period 3.
In this way we find the modular theta series
\[\sum_{\ell\in \Z^2}\{\sign(\ell_1)+\sign(\ell_2)\}\ko \left(\frac{-3}{\ell_1+\ell_2}\right) \ko q^{\ell_1^2+5\ell_1\ell_2+\ell_2^2},\] 
which we identify as $4\ko\eta_3\eta_{21}$.
Note that this is also an example given in \cite{polishchuk}, where a similar construction yields a big number of modular forms associated to indefinite binary quadratic forms.
We can also transfer our construction to other quadratic forms, for instance $Q(v)=v_1^2+6v_1v_2+v_2^2$, so that we obtain the theta series
\[\sum_{\ell\in \Z^2}\{\sign(\ell_1)+\sign(\ell_2)\}\ko \left(\frac{-4}{\ell_1+\ell_2}\right) \ko q^{\ell_1^2+6\ell_1\ell_2+\ell_2^2},\] which equals $4\ko\eta_8\eta_{16}$.
%We obtain similar examples for any $Q(v)=v_1^2+tv_1v_2+v_2^2$ with $t\geq 3$ by defining periodic functions $m_t:\Z^2\longrightarrow\C$ with $m_t\circ g=m_t$.
\end{example}

So for signature $(1,1)$ we obtain examples (some of which are already known) by employing periodic functions and the homogeneous polynomial $f\equiv 1$.
However, we will now have a look at quadratic forms of signature $(2,1)$, and here we will take %homogeneous and
(spherical) polynomials of higher degree to obtain a new interesting set of examples.
We'll first give two examples, where we can immediately apply Theorem \ref{theo2} and Corollary \ref{cor1}.
In the first of those we will basically determine all pairs of homogeneous polynomials satisfying the condition of Theorem \ref{theo2}, hence providing many cases where we obtain almost holomorphic cusp forms of weight $3/2+d$.
To make this more precise we consider the following lemma.

\begin{lemma}\label{lemiso}
Let $Q:\R^2\longrightarrow \R$ be a quadratic form of signature $(1,1)$ and let $g\in \operatorname{Aut}(Q,\R^2)$ with $\det g=1$ and $g \neq \pm I$.
Let $d\in\N$ be odd and let $U_d \subset \R[x_1,x_2]$ be the vector space of homogeneous polynomials of degree $d$.
Further, let $\Psi_g$ be the endomorphism of $U_d$ given by $\Psi_g(f):=f-f\circ g$.
Then $\Psi_g$ is an automorphism of $U_d$.
\end{lemma}

\begin{proof}
Over $\R$ we can split $Q$ as the product of two linear factors: $Q=h_1\cdot h_2$, where this decomposition is unique up to the order of the factors and multiplication by a scalar.
Since $Q\circ g =Q$ we thus have
\[ \begin{cases} h_1\circ g &= \lambda\ko h_2,\\
h_2 \circ g &= \lambda^{-1} h_1,\end{cases} \qquad \text{or} \qquad \begin{cases} h_1\circ g &= \lambda\ko h_1,\\
h_2 \circ g &= \lambda^{-1} h_2,\end{cases}\]
with $\lambda\in\R^*$.
The first situation cannot occur, since then the matrix of the linear map $U_1\longrightarrow U_1$, $h\mapsto h\circ g$ with respect to the basis $\{h_1,h_2\}$ of $U_1$ would be $\left(\begin{smallmatrix} 0 & \lambda^{-1}\\ \lambda&0\end{smallmatrix}\right)$, which has determinant $-1$.
However, with respect to the canonical basis the matrix is $g^\trans$, which has determinant 1.
Further, in the second situation we have $\lambda\neq \pm 1$, since $\lambda=\pm 1$ would imply $g =\pm I$, which we have excluded.
So we have
\[ h_1\circ g = \lambda\ko h_1\qquad \text{and}\qquad h_2 \circ g = \lambda^{-1} h_2,\]
with $\lambda\in\R^*$ and $\lambda\neq \pm 1$.
With a suitable change of variables we can write any $f\in U_d$ as a homogeneous polynomial of degree $d$ in $h_1$ and $h_2$, that is as a linear combination of the monomials $h_1^{\alpha_1} h_2^{\alpha_2}$ with $\alpha_1+\alpha_2=d$.
For such monomials we then have
\[ \Psi_g(h_1^{\alpha_1} h_2^{\alpha_2}) = (1-\lambda^{\alpha_1-\alpha_2})\ko h_1^{\alpha_1} h_2^{\alpha_2},\]
where $1-\lambda^{\alpha_1-\alpha_2}\neq 0$, since $\lambda \neq \pm 1$ and $\alpha_1\neq \alpha_2$ ($\alpha_1+\alpha_2=d$ is odd).
It now follows directly that $\Psi_g$ is an isomorphism.
For the inverse we have
\[ \Psi_g^{-1} (h_1^{\alpha_1} h_2^{\alpha_2}) = \frac1{1-\lambda^{\alpha_1-\alpha_2}}\ko  h_1^{\alpha_1} h_2^{\alpha_2}.\qedhere\]
\end{proof}

\begin{example}\label{ex284}
Let $Q(v)=v_1^2+4v_2^2-2v_3^2$ and fix $c=\frac{1}{\sqrt{2}}\left(\begin{smallmatrix}0\\0\\-1\end{smallmatrix}\right)\in \Ce_Q$.
Further, we pick the matrices
\[g_1=\begin{pmatrix} 1&0&0\\ 0&3&2\\0&4&3 \end{pmatrix} \qquad \text{and} \qquad g_2=\begin{pmatrix} 3&0&4\\ 0&1&0\\2&0&3 \end{pmatrix} \]
from the automorphism group $\aut{3}$.
For arbitrary odd degree $d$ we construct homogeneous polynomials $f_1$ and $f_2$ such that the condition $\sum_{i=1}^2 (f_i-f_i\circ g_i)=0$ from Theorem \ref{theo2} is fulfilled.
We start by noting that we can directly eliminate certain polynomials which will necessarily give vanishing theta series:
if $f_1$ is odd in the first variable $v_1$, then we have
\[\{\sign(v_3)-\sign(-4v_2+3v_3)\}\ko f_1(v)=-\{\sign(v_3)-\sign(-4v_2+3v_3)\}\ko f_1(v)\]
under the substitution $v_1\rightarrow -v_1$.
Therefore, we assume that $f_1$ is even in $v_1$ and hence so is $f_1-f_1\circ g_1$.
Similarly, we assume that $f_2$ and $f_2-f_2\circ g_2$ are even in $v_2$.
Thus we are looking for polynomials $f_1$ and $f_2$ for which $f_1-f_1\circ g_1=-(f_2-f_2\circ g_2)$ is in the vector space $V_d$ of homogeneous polynomials of degree $d$ that are even in both $v_1$ and $v_2$ (and odd in $v_3$).
We note that $\dim V_d=(d+1)(d+3)/8$.

The idea now is that we can start with any polynomial $f\in V_d$ and use Lemma \ref{lemiso}
to construct a unique polynomial $f_1$, which is even in $v_1$ and satisfies $f_1-f_1\circ g_1=f$:
we write $f(v)$ as $\sum_k v_1^{2k} p_k(v_2,v_3)$, where $p_k$ is a homogeneous polynomial of degree $d-2k$.
Then $\sum_k v_1^{2k} (\Psi_{\widetilde g_1}^{-1} (p_k))(v_2,v_3)$, where $\widetilde g_1=\left(\begin{smallmatrix} 3&2\\4&3\end{smallmatrix}\right)$, satisfies $f_1-f_1\circ g_1=f$.
Similarly, there exists a unique polynomial $f_2$, which is even in $v_2$ and satisfies $f_2-f_2\circ g_2=-f$.
This way we obtain a vector space of dimension $(d+1)(d+3)/8$ of solutions $(f_1,f_2)$ satisfying the condition $\sum_{i=1}^2 (f_i-f_i\circ g_i)=0$ from Theorem \ref{theo2}.

We now go a step further and assume that $f$ is spherical and consider the corresponding polynomials $f_1$ and $f_2$.
Since $f_1$ is even in $v_1$, so is $\Delta f_1$.
Further, we have $\Delta(f_1\circ g_1)= (\Delta f_1)\circ g_1$ and so we get
\[ (\Delta f_1) -(\Delta f_1)\circ g_1= \Delta( f_1-f_1\circ g_1)=\Delta f=0,\]
which by Lemma \ref{lemiso} and our previous construction yields $\Delta f_1=0$.
Similarly, we also have $\Delta f_2=0$.
Hence we have shown that if $f$ is spherical, then so are $f_1$ and $f_2$.
Further, the map $\Delta|_{V_d}:V_d\longrightarrow V_{d-2}$ is surjective, which one can easily check by using that $\Delta (v_1^{\alpha_1}v_2^{\alpha_2} v_3^{\alpha_3})$ is a linear combination of $v_1^{\alpha_1-2}v_2^{\alpha_2} v_3^{\alpha_3}$, $v_1^{\alpha_1}v_2^{\alpha_2-2} v_3^{\alpha_3}$ and $v_1^{\alpha_1}v_2^{\alpha_2} v_3^{\alpha_3-2}$, together with induction on $\alpha_1+\alpha_2$.
Therefore, the kernel of $\Delta|_{V_d}$ has dimension $\dim V_d - \dim V_{d-2}= (d+1)/2$, and so we obtain a vector space of dimension $(d+1)/2$ of spherical solutions $(f_1,f_2)$ satisfying the condition $\sum_{i=1}^2 (f_i-f_i\circ g_i)=0$ from Corollary \ref{cor1}. 

For $d=1$ the vector space $V_1$ is one-dimensional and is spanned by $f(v)=2v_3$.
The corresponding polynomials $f_1$ and $f_2$ are easily determined to be $f_1(v)=-2v_2+v_3$ and $f_2(v)=v_1-v_3$.
Since they are spherical, we obtain the following holomorphic cusp form of weight $5/2$ on $\Gamma_0(16)$:
\begin{align*}
\begin{split}
\sum_{i=1}^2&\modtheta^{c,g_i c}[f_i](\tau)=\sum_{i=1}^2\Theta^{c,g_i c}[f_i](\tau)\\
&=\sum\limits_{\ell\in\Z^3} \{\sign\bigl(B(c,\ell)\bigr)-\sign\bigl(B(g_1c,\ell)\bigr)\}\ko f_1(\ell)\ko q^{Q(\ell)}+\sum\limits_{\ell\in\Z^3} \{\sign\bigl(B(c,\ell)\bigr)-\sign\bigl(B(g_2c,\ell)\bigr)\}\ko f_2(\ell)\ko q^{Q(\ell)}\\
&=\sum\limits_{\ell\in\Z^3} \bigl\{\bigl(\sign(\ell_3)-\sign(-4\ell_2+3\ell_3)\bigr)\ko (-2\ell_2+\ell_3) + \bigl(\sign(\ell_3)-\sign(-2\ell_1+3\ell_3)\bigr)\ko (\ell_1-\ell_3)\bigr\}\ko q^{Q(\ell)}
\end{split}
\end{align*}
We identify this theta function as the eta product $4\ko \eta_2^2\eta_4\eta_8^2$.

For $d=3$ the vector space $V_3$ has dimension three.
In Table \ref{tabv3} we list a possible basis of polynomials $f$, together with $\Delta f$ and the corresponding polynomials $f_1$ and $f_2$.

\begin{table}[hbtp]
\renewcommand{\arraystretch}{1.4}
{\begin{tabular}{c|c|c|c}
$f$ & $\Delta f$ & $f_1$ & $f_2$  \\ \hline
$v_1^2 v_3$ & $v_3$ & $\frac12 v_1^2(-2v_2+v_3)$ & $\frac1{14}v_1(3v_1-4v_3)(v_1-v_3)$  \\ \hline
$v_2^2v_3$ & $\frac14 v_3$ & $\frac1{14}v_2(3v_2-2v_3)(-2v_2+v_3)$ & $\frac12 v_2^2(v_1-v_3)$  \\ \hline
$v_3^3$ & $-\frac32 v_3$ & $-\frac1{14}(8v_2^2+4v_2v_3-7v_3^2)(-2v_2+v_3)$ & $-\frac1{14}(2v_1^2+2v_1v_3-7v_3^2)(v_1-v_3)$ 
\end{tabular}}
\vspace{1em}
\caption{Basis elements $f$ of $V_3$ and the corresponding polynomials $f_1$ and $f_2$}
\label{tabv3}
\end{table}

In all three cases the corresponding theta function $\sum_{i=1}^2\modtheta^{c,g_i c}[f_i]=\sum_{i=1}^2\widehat\Theta^{c,g_i c}[f_i]$ is an almost holomorphic cusp form and the holomorphic theta function $\sum_{i=1}^2\Theta^{c,g_i c}[f_i]$ is a quasimodular form of weight 9/2 and depth 1 on $\Gamma_0(16)$.
We can identify these quasimodular forms as
\[\begin{split}
\frac87\ko \eta(2\tau)^2 \eta(4\tau) \eta(8\tau)^2 &\Bigl( G_2(\tau)-5 G_2(2\tau) + 10 G_2(8\tau) - 24 G_2(16\tau) +4\ko \frac{\eta(8\tau)^2\eta(16\tau)^4}{\eta(4\tau)^2}\Bigr),\\
-\frac27\ko \eta(2\tau)^2 \eta(4\tau) \eta(8\tau)^2 &\Bigl( G_2(\tau)- G_2(2\tau) -2G_2(4\tau) + 26 G_2(8\tau) - 24 G_2(16\tau) +4\ko \frac{\eta(8\tau)^2\eta(16\tau)^4}{\eta(4\tau)^2}\Bigr)
\end{split}\]
and
\[ \frac{12}7\ko \eta(2\tau)^2 \eta(4\tau) \eta(8\tau)^2 \bigl( G_2(2\tau) + 3 G_2(4\tau) +4 G_2(8\tau) \bigr),\]
where $G_2(\tau):=-\frac1{24} + \sum_{n=1}^\infty \sigma_1(n)\ko q^n$ is the quasimodular Eisenstein series of weight 2 and depth 1.

From Table \ref{tabv3} we can directly construct two linearly independent spherical solutions:
for the spherical polynomial $f(v)=(v_1^2-4v_2^2) v_3$ we have
\[ f_1(v) = \frac1{14}(7v_1^2-12v_2^2+8v_2v_3)(-2v_2+v_3)\qquad \text{and}\qquad  f_2(v)= \frac1{14}(3v_1^2-4v_1v_3-28v_2^2)(v_1-v_3).\]
The corresponding holomorphic theta function $\sum_{i=1}^2\Theta^{c,g_i c}[f_i]$ is a modular form of weight 9/2 on $\Gamma_0(16)$ and equals
\[ \frac{16}7\ko \eta(2\tau)^2 \eta(4\tau) \eta(8\tau)^2 \Bigl( G_2(\tau)- 3 G_2(2\tau) -G_2(4\tau) + 18 G_2(8\tau) - 24 G_2(16\tau) +4\ko \frac{\eta(8\tau)^2\eta(16\tau)^4}{\eta(4\tau)^2}\Bigr).\]
For the spherical polynomial $f(v)=(3v_1^2+12v_2^2+4 v_3^2)v_3$ we have
\[\begin{split}
f_1(v)&= \frac1{14}(21v_1^2+4v_2^2-40v_2v_3+28v_3^2)(-2v_2+v_3),\\
f_2(v)&= \frac1{14} (v_1^2+84v_2^2+28v_3^2-20v_1v_3)(v_1-v_3).
\end{split}\]
The corresponding holomorphic theta function is again modular and equals
\[ -\frac{48}7\ko \eta(2\tau)^2 \eta(4\tau) \eta(8\tau)^2 \bigl( G_2(2\tau) -4 G_2(4\tau) + 4 G_2(8\tau) \bigr).\]
\end{example}

\begin{example}
Considering the twelve diagonalized quadratic forms of signature $(2,1)$ and level 24 and choosing spherical polynomials of degree 1, we already obtain eight different eta quotients of weight $5/2$. We list possible choices of $Q$, $g_i$ and $f_i$ and the corresponding eta quotients that $\sum_{i=1}^2\modtheta^{c,g_i c}[f_i]$ evaluate to in Table \ref{tab24} (we omit the quadratic forms which lead to the same eta quotients).
As in Example \ref{ex284} one could generate many (almost) holomorphic cusp forms of weight $3/2+d$ by constructing suitable homogeneous polynomials of higher degree $d$.

\begin{table}[hbtp]
{\begin{tabular}{  m{7em} | m{12.5em} | m{10em} | m{10em} }
$Q$ & $g_i\in\aut{3}$ & $f_i$ spherical of\newline degree 1 & $\frac14\sum_{i=1}^2\modtheta^{c,g_i c}[f_i]$\newline $=\frac14\sum_{i=1}^2\Theta^{c,g_i c}[f_i]$\\
\hline
$v_1^2+3v_2^2-2v_3^2$ & $g_1=\left(\begin{smallmatrix} 1&0&0\\ 0&5&4\\ 0&6&5 \end{smallmatrix}\right),\,g_2=\left(\begin{smallmatrix} 3&0&4\\ 0&1&0\\ 2&0&3  \end{smallmatrix}\right)$ & $f_1(v)=-3v_2+2v_3,$\newline $f_2(v)=2v_1-2v_3$ & $\eta _2^3\eta_4^2\eta_6\eta_{24}/(\eta_8\eta_{12})$\\
\hline
 $v_1^2+2v_2^2-3v_3^2$ & $g_1=\left(\begin{smallmatrix} 1&0&0\\ 0&5&6\\ 0&4&5 \end{smallmatrix}\right),\,g_2=\left(\begin{smallmatrix} 2&0&3\\ 0&1&0\\ 1&0&2  \end{smallmatrix}\right)$  &
$f_1(v)=-v_2+v_3,$\newline $f_2(v)=v_1-v_3$ &  $\eta\eta_6^9\eta_8^2/(\eta_2\eta_3^3\eta_{12}^3)$\\
\hline
$v_1^2+6v_2^2-2v_3^2$ & $g_1=\left(\begin{smallmatrix} 1&0&0\\ 0&2&1\\ 0&3&2 \end{smallmatrix}\right),\,g_2=\left(\begin{smallmatrix} 3&0&4\\ 0&1&0\\ 2&0&3  \end{smallmatrix}\right)$  &
$f_1(v)=-3v_2+v_3,$\newline $f_2(v)=v_1-v_3$ &  $\eta_2^2\eta_3\eta_4^3\eta_{12}/(\eta\eta_6)$\\
\hline
$v_1^2+2v_2^2-6v_3^2$ & $g_1=\left(\begin{smallmatrix} 1&0&0\\ 0&2&3\\ 0&1&2 \end{smallmatrix}\right),\,g_2=\left(\begin{smallmatrix} 5&0&12\\ 0&1&0\\ 2&0&5  \end{smallmatrix}\right)$  &
$f_1(v)=-2v_2+2v_3,$\newline $f_2(v)=v_1-2v_3$ &  $\eta^2\eta_8\eta_{12}^9/(\eta_4\eta_6^3\eta_{24}^3)$\\
\hline
$v_1^2+6v_2^2-3v_3^2$ & $g_1=\left(\begin{smallmatrix} 1&0&0\\ 0&3&2\\ 0&4&3 \end{smallmatrix}\right),\,g_2=\left(\begin{smallmatrix} 2&0&3\\ 0&1&0\\ 1&0&2  \end{smallmatrix}\right)$  &
$f_1(v)=-2v_2+v_3,$\newline $f_2(v)=v_1-v_3$ &  $\eta_2\eta_6^3\eta_8\eta_{12}^2/(\eta_4\eta_{24})$\\
\hline
$3v_1^2+6v_2^2-v_3^2$ & $g_1=\left(\begin{smallmatrix} 1&0&0\\ 0&5&2\\ 0&12&5 \end{smallmatrix}\right),\,g_2=\left(\begin{smallmatrix} 2&0&1\\ 0&1&0\\ 3&0&2  \end{smallmatrix}\right)$  &
$f_1(v)=-3v_2+v_3,$\newline $f_2(v)=3v_1-v_3$ &  $\eta_2^9\eta_3\eta_{24}^2/(\eta^3\eta_4^3\eta_6)$\\
\hline
$2v_1^2+6v_2^2-3v_3^2$ & $g_1=\left(\begin{smallmatrix} 1&0&0\\ 0&3&2\\ 0&4&3 \end{smallmatrix}\right),\,g_2=\left(\begin{smallmatrix}5&0&6\\ 0&1&0\\ 4&0&5  \end{smallmatrix}\right)$  &
$f_1(v)=-2v_2+v_3,$\newline $f_2(v)=v_1-v_3$ &  $\eta\eta_4\eta_6^2\eta_{12}^3/(\eta_2\eta_3)$\\
\hline
$3v_1^2+6v_2^2-2v_3^2$ & $g_1=\left(\begin{smallmatrix} 1&0&0\\ 0&2&1\\ 0&3&2 \end{smallmatrix}\right),\,g_2=\left(\begin{smallmatrix}5&0&4\\ 0&1&0\\ 6&0&5  \end{smallmatrix}\right)$  &
$f_1(v)=-6v_2+2v_3,$\newline $f_2(v)=3v_1-2v_3$ &  $\eta_3^2\eta_4^9\eta_{24}/(\eta_2^3\eta_8^3\eta_{12})$\\
\end{tabular}}
\vspace{1em}
\caption{Eight different cusp forms of weight $5/2$ on $\Gamma_0(24)$}
\label{tab24}
\end{table}
\end{example}

Just like for quadratic forms of signature $(1,1)$, we obtain further examples for quadratic forms of signature $(2,1)$ if we modify the polynomials in the theta series by introducing an additional periodic factor as described in Remark \ref{remexamples}.

\begin{example}\label{ex2N}
Let $Q(v)=v_1^2+v_2^2-Nv_3^2$, where $N\in \N$ is such that $2N$ is not a perfect square.
As matrices in the automorphism group $\aut{3}$ we choose
\[g_1=\begin{pmatrix} -1&0&0\\ 0&1&0\\0&0&1 \end{pmatrix},\qquad g_2=\begin{pmatrix} 1&0&0\\ 0&-1&0\\0&0&1 \end{pmatrix}\qquad\text{and}\qquad g_3=\begin{pmatrix} \frac{x+1}2&\frac{x-1}2&Ny\\ \frac{x-1}2&\frac{x+1}2&Ny\\y&y&x \end{pmatrix}, \]
where $(x,y)$ is an integer solution to the Pell equation $x^2-2Ny^2=1$.
We take the spherical polynomials $f_1(v)=\frac12 y v_2 +\frac{x-1}4v_3$, $f_2(v)=\frac12 y v_1 +\frac{x-1}4v_3$ and $f_3(v)=v_3$ and consider the periodic function 
\[m(v)=\Bigl(\frac{-4}{v_1}\Bigr)\ko \Bigl(\frac{-4}{v_2}\Bigr),\]
where the Dirichlet character $\bigl(\frac{-4}{\cdot}\bigr)$ is given by 
\[\Bigl(\frac{-4}{n}\Bigr)=\begin{cases}
\phantom{-}1 &\text{if}\ n\equiv\phantom{-} 1 \ (\smod{4}),\\
-1&\text{if}\ n\equiv -1 \ (\smod{4}),\\
\phantom{-}0&\text{otherwise.}
\end{cases}\]
We set $\widetilde{f}_i:=m\cdot f_i$ and observe that $m\circ g_1 =m\circ g_2 =-m$ and $m\circ g_3=m$.
Hence $\sum_{i=1}^3(\widetilde{f}_i-\widetilde{f}_i\circ g_i)=0$ is equivalent to $f_1+f_1\circ g_1+f_2+f_2\circ g_2+f_3-f_3\circ g_3=0$, which we can easily verify for the polynomials $f_1$, $f_2$ and $f_3$ above.
Thus, $\sum_{i=1}^3\modtheta^{c,g_i c}[m,f_i]=\sum_{i=1}^3\Theta^{c,g_i c}[m,f_i]$ is a cusp form of weight $5/2$.
For $c=\frac1{\sqrt{N}}\Bigl(\begin{smallmatrix}0\\0\\-1\end{smallmatrix}\Bigr)$ we have $g_1c=g_2c=c$, so for this choice of $c\in\Ce_Q$ the first two theta functions $\Theta^{c,g_1 c}[m,f_1]$ and $\Theta^{c,g_2 c}[m,f_2]$ vanish.
Hence we obtain the modular theta series
\[\begin{split}\sum_{i=1}^3\modtheta^{c,g_i c}[m,f_i]&=\Theta^{c,g_3 c}[m,f_3](\tau)\\
&=\sum\limits_{\ell\in\Z^3} \bigl\{\sign(\ell_3)+\sign(y\ell_1+y\ell_2-x\ell_3)\bigr\}\ko \Bigl(\frac{-4}{\ell_1}\Bigr)\ko \Bigl(\frac{-4}{\ell_2}\Bigr)\ko \ell_3\ko q^{\ell_1^2+\ell_2^2-N\ell_3^2}.\end{split}\]
For $N=1$ we take $(x,y)=(3,2)$ as our solution to the equation $x^2-2y^2=1$.
The corresponding theta series is $4\ko \eta_2^5\eta_8^2/\eta^2$.
For $N=3$ we take $(x,y)=(5,2)$ and identify the theta series as
\[ 8\ko \frac{\eta_8^4\eta_{24}^9}{\eta_4\eta_{12}^3\eta_{16} \eta_{48}^3} - 32\ko \eta_8\eta_{16}\eta_{48}^3.\]
For $N=4$ we take $(x,y)=(3,1)$ and the theta series equals $4\ko \eta_4^2 \eta_8\eta_{16}^2$.
For $N=6$, $x=7$ and $y=2$ we find $8\ko \eta_2\eta_8^4\eta_{12}^2/(\eta_4\eta_6)$.

For $N=1$ and $d=3$ we could for example take the spherical polynomials
\[\begin{split}
f_1(v) &= 7(v_2+v_3)^3,\\
f_2(v) &= 7(v_1+v_3)^3,\\
f_3(v) &=4v_1^3+4v_2^3+3(v_1+v_2)v_3^2-9v_1v_2(v_1+v_2).
\end{split}\]
Again, the condition $f_1+f_1\circ g_1+f_2+f_2\circ g_2+f_3-f_3\circ g_3=0$ is satisfied, so
\[ \sum\limits_{\ell\in\Z^3} \bigl\{\sign(\ell_3)+\sign(2\ell_1+2\ell_2-3\ell_3)\bigr\}\ko \Bigl(\frac{-4}{\ell_1}\Bigr)\ko \Bigl(\frac{-4}{\ell_2}\Bigr)\ko f_3(\ell)\ko q^{\ell_1^2+\ell_2^2-\ell_3^2}\]
is a cusp form of weight $9/2$, which we can identify as
\[ 48\ko \frac{\eta(2\tau)^5\eta(8\tau)^2}{\eta(\tau)^2}\bigl( G_2(\tau)-5\ko G_2(2\tau) +12\ko G_2(8\tau)\bigr).\]
\end{example}

In the next example, we give a similar construction for quadratic forms of level $4N$.
Here we assume that $N$ itself is not a perfect square, and consider $I=\{1,2\}$ and a different periodic function.
%leads to another set of examples.

\begin{example}\label{exN}
Let $Q(v)=v_1^2+v_2^2-Nv_3^2$, where $N\in \N$ is not a perfect square. As matrices in the automorphism group we choose \[g_1=\begin{pmatrix} 0&1&0\\ 1&0&0\\0&0&1 \end{pmatrix}\qquad\text{and}\qquad g_2=\begin{pmatrix} 1&0&0\\ 0&x&Ny\\0&y&x \end{pmatrix},\]
where $(x,y)$ denotes an integer solution of the Pell equation $x^2-Ny^2=1$ for which $x$ is odd (if for the fundamental solution $(x_1,y_1)$ the integer $x_1$ is even, we take the solution $(x_2,y_2)$ given by $x_2=x_1^2+Ny_1^2=1+2Ny_1^2$, $y_2=2x_1y_1$).
%This construction of $g_2$ also ensures that $g_2$ has infinite order.
Further, we choose $c=\frac{1}{\sqrt{N}}\Bigl(\begin{smallmatrix}0\\0\\-1\end{smallmatrix}\Bigr)\in \Ce_Q$. As $g_1c=c$ holds, the theta series $\modtheta^{c,g_1 c}$ vanishes, thus we only need to determine $\sign\bigl(B(c,v)\bigr)=\sign(v_3)$ and 
$\sign\bigl(B(g_2 c,v)\bigr)=-\sign(yv_2-xv_3)$ to write out the theta series later.
For the periodic function 
\begin{align*}
m(v)=\begin{cases}
(-1)^{v_2}&\text{if}\ v_1\not\equiv v_2 \ (\smod{2}),\\
0&\text{if}\ v_1\equiv v_2 \ (\smod{2})
\end{cases}
\end{align*}
we have $m\circ g_1=-m$ and $m\circ g_2=m$, as $x$ is odd and $y$ is even.
For a polynomial $f_i$, we again set $\widetilde{f}_i=m\cdot f_i$ and observe that $\sum_{i=1}^2 (\widetilde{f}_i-\widetilde{f}_i\circ g_i)=0$ is equivalent to $f_1+f_1\circ g_1+f_2-f_2\circ g_2=0$. 

Let $\alpha:=y/\operatorname{gcd}(x-1,y)$ and $\beta:=(x-1)/\operatorname{gcd}(x-1,y)$.
For $0\leq r\leq (d-1)/2$, we define the homogeneous polynomials of degree $d$
\[f_1(v):=f^r_1(v)=h^r_1(v)\ko Q(v)^r\qquad\text{and}\qquad f_2(v):=f^r_2(v)=h^r_2(v)\ko Q(v)^r,\] where
\[\begin{split}
h^r_1(v)&=\frac14\Big\{(\alpha(v_1+ v_2)+\beta v_3)^{d-2r}+(\alpha(-v_1+ v_2)+\beta v_3)^{d-2r}\\
&\qquad+(\alpha(v_1- v_2)+\beta v_3)^{d-2r}+(\alpha(-v_1- v_2)+\beta v_3)^{d-2r}\Big\},\\
h^r_2(v)&=\frac12\Big\{(\alpha(v_1+ v_2)-\beta v_3)^{d-2r}+(\alpha(-v_1+ v_2)-\beta v_3)^{d-2r}\Big\}.
\end{split}\]
By a simple calculation, where we use that $f_1=f_1\circ g_1$, $d$ is odd and $g_2$ leaves $Q$ invariant, we check that the polynomials above fulfill the equation $f_1+f_1\circ g_1+f_2-f_2\circ g_2=0$. We use these to construct the almost holomorphic cusp forms
\[\sum_{i=1}^2\modtheta^{c,g_ic}[m,f_i](\tau)=\widehat\Theta^{c,g_2 c}[m,f_2](\tau)=y^{-d/2}\sum\limits_{\ell\in\Z^3} \bigl\{\sign(\ell_3)+\sign(y\ell_2-x\ell_3)\bigr\}\ko m(\ell)\ko \widehat{f}_2(\ell y^{1/2})\ko q^{Q(\ell)}\]
of weight $3/2+d$. In general, $f_1$ and $f_2$ are not spherical, but applying the Laplacian, we see that for $i=1,2$
\[\Delta \bigl(f^r_i(v)\bigr)=r\ko(2d-2r+1)\ko h^r_i(v)\ko Q(v)^{r-1}
+\frac12 (2\alpha^2-\beta^2/N)\ko(d-2r)\ko(d-2r-1)\ko h^{r+1}_i(v)\ko Q(v)^r\]
holds. So we can choose a linear combination of these homogeneous polynomials to construct the spherical polynomials
%\[P_i:=\binom{2d}{d}^{-1}\cdot\sum\limits_{r=0}^{(d-1)/2} \bigl(\beta^2/N-2\alpha^2\bigr)^r\ko \binom{2(d-r)}{r,d-r,d-2r}\ko p^r_i\quad(i=1,2).\]
\[F_i:=\sum\limits_{r=0}^{(d-1)/2} \bigl(\beta^2/N-2\alpha^2\bigr)^r\ko \binom{2(d-r)}{r,d-r,d-2r}\ko f^r_i\quad(i=1,2),\]
using the notation of trinomial coefficients $\binom{n}{j,k,\ell}:=\frac{n!}{j!\ko k!\ko \ell!}$ with $n=j+k+\ell$ for $j,k,\ell\in \N_0$.
Hence, we obtain the holomorphic cusp form
\[\sum_{i=1}^2\modtheta^{c,g_ic}[m,F_i](\tau)=\Theta^{c,g_2 c}[m,F_2](\tau)=\sum\limits_{\ell\in\Z^3} \bigl\{\sign(\ell_3)+\sign(y\ell_2-x\ell_3)\bigr\}\ko m(\ell)\ko F_2(\ell )\ko q^{Q(\ell)}.\]
For $N=2$ we take $(x,y)=(3,2)$, thus $\alpha=\beta=1$. Then we have for $d=1$ the spherical polynomials $f_1(v)=v_3$ and $f_2(v)=v_2-v_3$ and the theta series equals $(-4)\ko \eta_2^2\eta_4\eta_8^2$ (we obtained the same eta product in Example \ref{ex284}).

For $d=3$ we obtain for $r=0$ the pair of homogeneous polynomials
\[\begin{split}
f^0_1(v)&=v_3\bigl(3v_1^2+3v_2^2+v_3^2\bigr),\\
f^0_2(v)&%=3v_1^2v_2+v_2^3-3v_1^2v_3-3v_2^2v_3+3v_2v_3^2-v_3^3
= 3v_1^2(v_2-v_3)+(v_2-v_3)^3.
\end{split}\]
The corresponding almost holomorphic theta series can be identified as
\[24\ko \eta(2\tau)^2\eta(4\tau)\eta(8\tau)^2\bigl(G_2^\ast(2\tau)-G_2^\ast(4\tau)+4G_2^\ast(8\tau)\bigr),\]
where $G_2^*$ is the almost holomorphic modular form of weight 2 given by $G_2^\ast(\tau):=G_2(\tau)+\frac1{8\pi y}$.
For $r=1$ we have 
\[\begin{split}
f^1_1(v)&=v_3\bigl(v_1^2+v_2^2-2v_3^2\bigr),\\
f^1_2(v)&%=v_1^2v_2+v_2^3-v_1^2v_3-v_2^2v_3-2v_2v_3^2+2v_3^3
=(v_2-v_3)(v_1^2+v_2^2-2v_3^2)
\end{split}\]
and the corresponding theta series is
\[16\ko \eta(2\tau)^2\eta(4\tau)\eta(8\tau)^2\bigl(G_2^\ast(2\tau)+G_2^\ast(4\tau)+4G_2^\ast(8\tau)\bigr).\]
Thus, we have constructed two almost holomorphic cusp forms of weight $9/2$ and depth 1. %level 16

Now we set $F_i=10f^0_i-9f^1_i$ (we take the formula for $F_i$ given above, but divide by the greatest common divisor of the integer coefficients in the linear combination) and get
the holomorphic cusp form
\[96\ko \eta(2\tau)^2\eta(4\tau)\eta(8\tau)^2\bigl(G_2(2\tau)-4G_2(4\tau)+4G_2(8\tau)\bigr).\]

For $N=3$ we take $(x,y)=(7,4)$, thus $\alpha=2$ and $\beta=3$. For $d=1$ we have the spherical polynomials $f_1(v)=3v_3$ and $f_2(v)=2v_2-3v_3$ and the theta series equals $(-8)\ko \eta \eta_4^4\eta_6^2/(\eta_2\eta_3)$.

For $d=3$ we obtain for $r=0$ the pair of homogeneous polynomials
\[\begin{split}
f^0_1(v)&=9v_3\bigl(4v_1^2+4v_2^2+3v_3^2\bigr),\\
f^0_2(v)&%=24v_1^2v_2+8v_2^3-36v_1^2v_3-36v_2^2v_3+54v_2v_3^2-27v_3^3.
=12v_1^2(2v_2-3v_3)+(2v_2-3v_3)^3.
\end{split}\]
The corresponding almost holomorphic theta series can be identified as
\[24\ko \frac{\eta(\tau)\eta(4\tau)^4\eta(6\tau)^2}{\eta(2\tau)\eta(3\tau)}\bigl(3G_2^\ast(\tau)-6G_2^\ast(2\tau)-9G_2^\ast(3\tau)+8G_2^\ast(4\tau)+36G_2^\ast(6\tau)\bigr).\]
For $r=1$ we have 
\[\begin{split}
f^1_1(v)&=3v_3\bigl(v_1^2+v_2^2-3v_3^2\bigr),\\
f^1_2(v)&%=2v_1^2v_2+2v_2^3-3v_1^2v_3-3v_2^2v_3-6v_2v_3^2+9v_3^3
=(2v_2-3v_3)\bigl(v_1^2+v_2^2-3v_3^2\bigr)
\end{split}\]
and the corresponding theta series is
\[8\ko \frac{\eta(\tau)\eta(4\tau)^4\eta(6\tau)^2}{\eta(2\tau)\eta(3\tau)}\bigl(G_2^\ast(\tau)-2G_2^\ast(2\tau)-3G_2^\ast(3\tau)+16G_2^\ast(4\tau)+12G_2^\ast(6\tau)\bigr).\]
We set $F_i=f^0_i-3f^1_i$ and get
the holomorphic cusp form of weight $9/2$ %and level 48
\[48\ko \frac{\eta(\tau)\eta(4\tau)^4\eta(6\tau)^2}{\eta(2\tau)\eta(3\tau)}\bigl(G_2(\tau)-2G_2(2\tau)-3G_2(3\tau)-4G_2(4\tau)+12G_2(6\tau)\bigr).\]
\end{example}

%\section*{Acknowledgements}

\end{document}